\documentclass{amsart}

\usepackage{amsthm,amsfonts,amsmath,amssymb,latexsym,enumitem}
\usepackage{upref,amssymb,eucal,ae}
\usepackage[all,cmtip]{xy}

\newtheorem{theorem}{Theorem}

\newtheorem{corollary}[theorem]{Corollary}

\newenvironment{remark}{\medskip \refstepcounter{theorem}
\noindent  {\bf Remark \thetheorem}.\rm}{\,}

\makeatletter
\def\Ddots{\mathinner{\mkern1mu\raise\p@
\vbox{\kern7\p@\hbox{.}}\mkern2mu
\raise4\p@\hbox{.}\mkern2mu\raise7\p@\hbox{.}\mkern1mu}}
\makeatother

\newcounter{spthe}

\def\ddb{\partial  \overline{\partial}}

\def\<{\langle}
\def\a{\alpha}
\def\>{\rangle}

\def\tm{\tilde{M}}
\def\tn{\tilde{n}}
\def\tg{\tilde{g}}
\def\mb#1{{\mathbb #1}}
\def\mc#1{{\mathcal #1}}

\def\BOne{{\mathchoice {\rm 1\mskip-4mu l} {\rm 1\mskip-4mu l}
                          {\rm 1\mskip-4.5mu l} {\rm 1\mskip-5mu l}}}

\begin{document}
\title[Invariants of isometric embeddings of metrics on a surface]{Conformal
invariants of isometric embeddings of the smooth metrics on a surface}
\author{Santiago R. Simanca}
\dedicatory{To Wojciech Kucharz}
\thanks{Most of this research was carried out while the 
author was supported by the Simons Foundation Visiting Professorship award 
657746 at CIMS}
\email{srsimanca@gmail.com}

\begin{abstract}
We view all smooth metrics $g$ on a closed surface $\Sigma$ through their Nash 
isometric embeddings $f_g: (\Sigma,g) \rightarrow (\mb{S}^{\tn},\tg)$ into a 
standard sphere of large, but fixed, dimension $\tn$, and analyze the extrinsic
quantities of these embeddings corresponding to variations of the metrics 
within a fixed conformal class. We define the Willmore functional 
$\mc{W}_{f_g}$ over this space of metrics on $\Sigma$ in terms of the extrinsic 
quantities of $f_g$. Its infimum over metrics in a conformal class is an 
invariant of the class varying differentiably with it. We identify 
geometrically one 
class for which this invariant is the smallest, set the scale using its 
Einstein representative of minimal isometric embedding, and compute the 
invariant of an arbitrary class by using its Einstein representative of area 
equal to this scale. If $\Sigma$ is oriented 
of genus $k$, when $k=0$, we use the gap theorem of Simons to show that there 
is a unique conformal class of metrics on $\Sigma$, whose invariant $16\pi$
is the value for the standard totally geodesic 
embedding of $\mb{S}^2 \hookrightarrow \mb{S}^{\tn}$, and we have that
$\mc{W}_{f_g}(\Sigma) \geq 16\pi$, with the lower bound achieved if, and only
if, $f_g$ is conformally equivalent to this standard geodesic embedding 
of $\mb{S}^2$, and ${\rm area}_g(\Sigma) \leq 4\pi$,  
while when $k\geq 1$, the Lawson 
minimal surface $(\xi_{k,1},g_{\xi_{k,1}})$ fixes the scale, and we show that 
$\mc{W}_{f_g}(\Sigma) \geq 4\, {\rm area}_{g_{\xi_{k,1}}} (\xi_{k,1})$, with 
the lower bound achieved by $f_g$ if, and only if, $f_g$ is 
conformally equivalent to 
$f_{g_{\xi_{k,1}}} : (\xi_{k,1}, g_{\xi_{k,1}}) \rightarrow (\mb{S}^3,\tg)
\hookrightarrow (\mb{S}^{\tn},\tg)$, and ${\rm area}_g(\Sigma)\leq
{\rm area}_{g_{\xi_{k,1}}} (\xi_{k,1})$. For a nonoriented $\Sigma$, we prove a 
likewise estimate from below for $\mc{W}_{f_g}(\Sigma)$, and characterize 
conformally the surface that realizes the optimal lower bound.
\end{abstract}

\maketitle

\section{Introduction: On the Willmore conjecture} \label{s1}
Let $\Sigma $ be a closed surface, and 
$$ 
\iota : \Sigma \rightarrow S_c^{\tn} 
$$ 
be an immersion of $\Sigma$ into the $\tn$-dimensional simply connected 
space form $S_c^{\tn}$ of curvature $c$. The immersion confers to $\Sigma$ a
particular conformal structure induced by the metric on the background, which
can thus be associated with it.  If $H_{\iota}$ and 
$K^{\tilde{g}}(e_i,e_j)$ denote the mean curvature vector of $\iota$, 
and sectional curvature of $S_c^{\tn}$ spanned by $\{e_i,e_j\}$ at points in 
$\iota(\Sigma)$, respectively, we have the functionals  
\begin{eqnarray} 
\Psi_{\iota}(\Sigma) & = & 
{\displaystyle \int_{\iota(\Sigma)} \| H_{\iota} \|^2 d\mu_{\iota(\Sigma)}}\, , 
\vspace{1mm} \label{fu1} \\ \Theta_{\iota}(\Sigma) & = 
& {\displaystyle \int_{\iota(\Sigma)} \sum_{i,j}K^{\tg}(e_i,e_j) d\mu_{\iota(
\Sigma)}} \label{fu2}\, , 
\end{eqnarray}
where $d\mu_{\iota(\Sigma)}$ is the area measure on $\iota(\Sigma)$. The 
Willmore energy of $(\Sigma, \iota)$ is then defined by 
\begin{equation} \label{we}
\mc{W}_{\iota}(\Sigma)= 2\Theta_{\iota}(\Sigma)+
\Psi_{\iota}(\Sigma) = 
\int_{\iota(\Sigma)} ( 4c + \| H_{\iota}\|^2) 
d\mu_{{\iota}(\Sigma)} \, .
\end{equation}

As discovered by Blasche using $\mb{R}^3$ as the background \cite{bla}, the 
Willmore energy is invariant under conformal transformations of the  
space form where it is defined, no matter its dimension. Since Euclidean 
space, hyperbolic space, and the punctured sphere are conformally equivalent, 
we may adjust the background harmlessly by a dilation, and reduce additional 
considerations to the cases $c=-1$, $c=0$, or $c=1$, respectively, or the
latter specifically given its advantages in the analysis of
the extrinsic quantities of immersions, and for which the density of 
$\mc{W}_{\iota}(\Sigma)$ is pointwise bounded below by the constant $4$.  

If $\Sigma=\Sigma^k$ is a a surface of topological genus $k$, we define the set 
$$
\mc{I}_{\mb{S}^3}(\Sigma^k)=\{\iota: \; \iota : \Sigma^k \rightarrow \mb{S}^3 
\text{\mbox{} is an immersion}\}\, ,  
$$
of all possible immersions of $\Sigma$ into $\mb{S}^3$. 
Other than for $\Sigma =\mb{P}^2(\mb{R})$, $\mc{I}_{\mb{S}^3}(\Sigma^k)$ is 
nonempty \cite[Corollary 1.6, Theorem 5]{la2}, and has distinguished elements 
in it that are minimal. If $\Sigma$ is oriented of genus $k=0$, and minimally
immersed into $\mb{S}^3$, the immersion is totally geodesic \cite{alm},
and $\iota(\Sigma) = \mb{S}^2 \hookrightarrow \mb{S}^3$ is the standard sphere.
If the genus $k\geq 1$, there exists an embedding 
$\xi_{k,1}=\iota_{\xi_{k,1}}(\Sigma)\hookrightarrow \mb{S}^3$ that
is minimal, and of equilateral fundamental domain \cite[Proposition 6.1]{la2}.
(We shall refer to $\xi_{k,1}$ as the equilateral Lawson Riemann surface of 
genus $k$, and denote by $g_{\xi_{k,1}}$ its intrinsic Riemannian metric.)
When $\Sigma$ is not orientable of genus $k>2$, there is also a {\it similarly}
distinguished minimal surface $\eta_{k-1,1}=\iota_{\eta_{k-1,1}}(\Sigma) 
\hookrightarrow \mb{S}^3$ of Lawson, but in general, it has  
self-intersections \cite[\S 8]{la2}. 
$\mb{P}^2(\mb{R})$ can be minimally embedded into $\mb{S}^4$ by a map from
the sphere $\mb{S}^2(r)$ of radius $r=\sqrt{3/2}$ that is a   
$2$-to-$1$ cover of the embedding, with the conformal structures 
of base and cover consistent with each other \cite{cdck,si3}.

We observe that $\xi_{1,1}=\mb{S}^2(1/\sqrt{2})\times \mb{S}^2(1/\sqrt{2})$
is the Clifford torus with its linear  embedding into $\mb{S}^3$, and that
$\mu_{g_{\xi_{1,1}}}(\xi_{1,1}) = 2\pi^2$.  
There is no explicit expression for $\mu_{g_{\xi_{k,1}}}(\xi_{k,1})$ if $k>1$,
but $\mu_{g_{\xi_{k,1}}}(\xi_{k,1})< 8\pi$, and 
$\mu_{g_{\xi_{k,1}}}(\xi_{k,1})\nearrow 8\pi$ as
$k \nearrow \infty$ \cite{ku,kls}.

For $\Sigma$ other than $\mb{P}^2(\mb{R})$, we define the 
Willmore functional over this domain, 
\begin{equation} \label{fun}
\begin{array}{ccl}
\mc{W} : \mc{I}_{\mb{S}^3}(\Sigma^k) & \rightarrow & \mb{R} \\
\mbox{}\hspace{6mm} \iota & \mapsto & \mc{W}_{\iota}(\Sigma^k) 
\end{array} \, .
\end{equation}
Assuming that $\Sigma$ is oriented, we recall a conjecture, 
posed originally by Willmore when $k=1$ \cite{wi}, and that we state as follows:
\medskip

\noindent {\bf Generalized Willmore Conjecture} (GWC). {\it If $\Sigma^k$ is an 
oriented surface of topological genus $k\geq 1$, then 
$$
\mc{W}_{\iota}(\Sigma^k) \geq 4 \mu_{g_{\xi_{k,1}}}
(\xi_{k,1}) \, , 
$$
and if $(\iota,\Sigma)$ achieves the lower bound, $\iota(\Sigma)$ is 
conformally equivalent to $(\xi_{k,1},g_{k,1})$ in $\mb{S}^3$.} 
\smallskip

The case of $k=1$ in this statement enhances the original conjecture, 
formulated using the conformally equivalent background $\mb{R}^3$ instead, by 
characterizing the conformal class of the model $\Sigma^{k=1}$ that achieves 
the lower bound, a later addendum of White \cite{whi}. The 
existence of a surface minimizer for any $k \geq 1$ was proved by Simon
\cite{lsim}, and Bauer and Kuwert \cite{baku}. When $k>1$, the optimal lower 
bound above is believed to be true commonly, though we have seen it so 
written only by Simon himself in \cite{lsim}, where he refers to it as a 
``tempting conjecture'' to make.

The original conjecture was solved by Cod\'a and Neves recently 
\cite{cone}, who established a stronger result clarifying 
an important aspect of the GWC in passing. Indeed, after an 
earlier contribution to the problem by Li and Yau \cite{liya}, which implies 
that a compact surface with self-intersections has Willmore energy that is at 
least $32 \pi$, it was realized that in the analysis of the optimal lower bound 
for $\mc{W}$, it is enough to consider its values over embedded surfaces only, 
as opposed to the more general immersed ones. Cod\'a and Neves proved that if 
$\Sigma $ is any embedded closed surface in $\mb{S}^3$ of genus $k\geq 1$, 
hence oriented,
then $\mc{W}_{\iota}(\Sigma) \geq 8\pi^2$, and that the equality holds if, 
and only 
if, $\Sigma$ is the Clifford torus up to conformal mappings of the background
\cite[Theorem A]{cone}. Thus, while proving the $k=1$ conjecture, they find 
also a lower bound for the Willmore energy of embedded surfaces of arbitrary 
genus $k\geq 1$. But elaborating on their proof to elucidate the conjecture in 
its general form when $k>1$ seems difficult.  

Cod\'a and Neves' proof is based on the min-max theory for the area functional 
of Almgren \cite{alm1}, and the implications of the fundamental min-max theorem 
of Pitts \cite{pitts}, which ties this theory with that of minimal surfaces in 
a crucial way. For let $\mc{Z}_2(\mb{S}^3)$ be the space of integral two 
currents in $\mb{S}^3$ that have trivial boundary, conveniently topologized. 
If $\Phi: [0,1]^k \rightarrow \mc{Z}_2(\mb{S}^3)$ is 
a continuous mapping, we denote by $[\Phi]$ the set of all continuous mappings 
$\tilde{\Phi}: [0,1]^k \rightarrow \mc{Z}_2(\mb{S}^3)$ homotopic to $\Phi$ 
through homotopies that are fixed on $\partial [0,1]^k$, and define
$$
L([\Phi])= \inf_{\tilde{\Phi}\in [\Phi]}\sup_{x\in [0,1]^k} {\rm area}
(\tilde{\Phi}(x)) \, .
$$
If 
$$
L([\Phi])> \sup_{x\in \partial [0,1]^k} {\rm area} (\tilde{\Phi}(x)) \, ,
$$
Pitts' theorem asserts the existence of a disjoint collection of smooth
surfaces $\Sigma_1, \ldots, \Sigma_N$, minimally embedded into $\mb{S}^3$,  
such that
$$
L([\Phi])= \sum_{i=1}^N m_i \, {\rm area}(\Sigma_i)
$$
for some positive integers $m_1, \ldots, m_N$. 
Given a surface $\Sigma$ of positive genus embedded in $\mb{S}^3$, 
Cod\'a and Neves produce an associated ``canonical'' path of
conformal mappings defined for $(x,t)\in [0,1]^4 \times [0,1]=[0,1]^5$, and 
apply Pitt's result to its homotopy class. The Clifford torus appears as the 
minimizer of the $L$ min-max function of this class, with the lower bound 
equal to that conjectured by Willmore. Notwithstanding a number of technical 
details that are worked out with dexterity, the hint that they should thus 
obtain the Clifford torus came from a result of Urbano \cite{urb}, who 
had characterized the surfaces involved in Simons' gap theorem 
\cite[Theorem 5.3.2, Corollary 5.3.2]{si} \cite[Main Theorem]{cdck} 
\cite[Corollary 2]{law0} at the upper and lower end as the only compact minimal 
surfaces of $\mb{S}^3$ of index no greater than five, with the equatorial 
sphere having index one, and the Clifford torus having index five.

The density of the Willmore functional over oriented 
surfaces $\Sigma$ embedded into  
$S^{\tn}_c$ is $4c +\| H_{\iota}\|^2$, so its optimal lower bound over $\Sigma$s
of positive genus in $\mb{S}^3$, if
achieved by a minimal surface, must be achieved by one that in its intrinsic 
metric is an extremal of the area function across conformal classes. 
By the properties of the Lawson surfaces $\xi_{k,1}$, we would expect that a 
minimizer should be the Clifford torus $\xi_{1,1}$. It is less clear that 
the intrinsic metric on this minimizer is Einstein, and even after we realize 
that this is so, it is still unclear how the {\it scale} of this Einstein torus 
is fixed. We clarify these issues by being far less ambitious 
in scope than Cod\'a and Neves in the analysis of $\mc{W}$, 
first looking at its range when the functional is defined over the space of all 
embeddings of a surface of fixed topology only, the genus that is, 
and then tying the minimality 
of the minimizer of $\mc{W}$ to an intrinsic conformal geometric property of 
canonical metric representatives of the class that this minimizer represents. 
We outline our approach, and state our main result next. 

We begin by recalling that if $(M^n,g)$ is a closed Riemannian manifold, 
and 
\begin{equation} \label{emb}
f_g: (M,g) \rightarrow (\tm^{\tn}, \tg)
\end{equation}
is an isometric embedding of it into some Riemannian background 
$(\tm^{\tn}, \tg)$, the scalar curvature $s_g$ of $g$ relates to the
extrinsic quantities of the embedding by 
\begin{equation}
s_g  = \sum K^{\tg}(e_i,e_j) +\tg(H,H)- \tg(\a,\a) \, .
\label{sce}
\end{equation}
Here, $\{ e_i\}$ is an orthonormal base of the tangent space to the
submanifold $f_g(M)$, $K^{\tg}$ is the sectional curvature of the background
metric $\tg$, and $\tg(H,H)$ and $\tg(\a,\a)$ are the squared $\tg$-norms of
the mean curvature $H$, and second fundamental form $\a$ tensors of $f_g$, 
respectively. We consider the extrinsic functionals $\Psi_{f_g}(M)$ and 
$\Theta_{f_g}(M)$ of earlier in this context, and if in addition we 
introduce the functional
\begin{equation} \label{fu3}
 \Pi_{f_g}(M)  =  
{\displaystyle \int_{f_g(M)} \| \alpha \|^2 d\mu_g } 
\end{equation}
associated to this embedding, we have the intrinsic total scalar curvature 
given by 
\begin{equation} \label{tsc}
\mc{S}_g(M):=\int_M s_g d\mu_g =\Theta_{f_g}(M)+\Psi_{f_g}(M) -\Pi_{f_g}(M)\, .
\end{equation}

By the Nash embedding theorem \cite{nash}, any Riemannian manifold $(M^n,g)$ 
can be isometrically embedded into a space form of sufficiently large but 
fixed dimension $\tn$, which depends on $n$ but not on $M$, and so 
any $S_c^{\tn}$ can be cast in the role of $(\tm^{\tn},\tg)$ above;  
we conveniently use the positive case because then, in addition to their 
compactness and large conformal group, the possibility of having minimal
embeddings of $(M,g)$ is not a priori obstructed. Hence, we identify the 
open cone $\mc{M}(M)$ of smooth Riemannian metrics on $M$ with the space 
of their isometric
embedded images into $(\mb{S}^{\tn},\tg) \subset (\mb{R}^{\tn+1}, \tg)$,
\begin{equation} \label{met}
\mc{M}_{\mb{S}^{\tn}}(M)=\{f_g(M)\subset \mb{S}^{\tn}:\; f_g: (M,g) 
\rightarrow (\mb{S}^{\tn},\tg) \; \text{isometric embedding}\}\, .  
\end{equation}
By Palais' extension theorem \cite{pal}, different metrics 
correspond to submanifolds that are ambiently isotopic to each 
other, with a submanifold that is ambiently isotopic to itself being 
associated to the same metric, in which case the said isotopy reflects the 
action of the diffeomorphism group on it. Since $(\mb{S}^{\tn},\tg)$ 
carries isometric embeddings of all the metrics on $M$, including those 
preserving volume in the same conformal class,
if $f_g \in \mc{M}_{\mb{S}^{\tn}}(M)$ the function $\| H_{f_g}\|^2$ is 
constant \cite[Theorem 6]{scs}, and the embeddings of smooth deformations of 
metrics correspond to 
smooth isotopic submanifolds whose mean curvature vectors have constant 
squared norm.
If $f$ is an embedding of $M$ into $\mb{S}^{m} 
\subset \mb{S}^{\tn}$, and $g$ is the induced metric on $f(M)$, 
if $\| H_f\|^2$ is not constant, there are volume preserving metrics on $M$ 
arbitrarily close to $g$ in the $C^2$ topology, which do not admit 
isometric embeddings into $(\mb{S}^{\tn},\tg)$ that are isotopic deformations 
of $f$. 

If $\Sigma^k$ is a surface of topological genus $k$, we define its Willmore 
functional by
\begin{equation} \label{will}
\begin{array}{rcl}
\mc{W}: \mc{M}_{\mb{S}^{\tn}}(\Sigma^k) & \rightarrow & \mb{R} \\
f_g\hspace{4mm}\mbox{}  & \rightarrow & \mc{W}_{f_g}(\Sigma) =2\Theta_{f_g}(
\Sigma) + \Psi_{f_g}(\Sigma={\displaystyle
\int_{f_g(\Sigma)} (4 + \| H_{f_g}\|^2) d\mu_{f_g(\Sigma)}} 
\end{array} \, .
\end{equation}
With the $C^2$-topology in its domain, this function is continuous;   
it is differentiable in the $C^l$-topology for any $l\geq 3$.

If $F: (\mb{S}^{\tn},\tg) \rightarrow (\mb{S}^{\tn},\tg)$ is a conformal 
diffeomorphism   
and $f_g(\Sigma) \hookrightarrow (\mb{S}^{\tn},\tg)$, then 
$F(f_g(\Sigma)) \hookrightarrow (\mb{S}^{\tn}, \tg)$ has intrinsic metric 
$\tg\mid_{F(f_g(\Sigma))}$, and $\mc{W}_{f_g}(\Sigma)= 
\mc{W}_{F\circ f_{g}}(\Sigma)$. We define
\begin{equation} \label{min}
\mc{W}(\Sigma,[g])=\inf_{g'\in[g]}\mc{W}_{f_{g'}}(\Sigma) =
\inf_{g'\in[g]} \int_{f_{g'}(\Sigma)} \left( 4 + \| H_{f_{g'}}\|^2 \right)
d\mu_{f_{g'}} \, ,  
\end{equation}
There are $g$s in $[g]$ such that 
$\mc{W}_{f_g}(\Sigma) =\mc{W}(\Sigma,[g])$, and as a function of $[g]$,
$\mc{W}(\Sigma,[g])$ varies differentiably. We find the $[g]$ that 
minimizes it. 

Critical points of the Willmore functional 
and normal stationary critical points of the 
functionals $\Psi_{f}(\Sigma)$ and $\Theta_f(\Sigma)$ in (\ref{fu1}) and
(\ref{fu2}), respectively, are closely related. The latter are 
isometric embeddings $f: (\Sigma,g) 
\rightarrow (\tm^{2+q}, \tg)$ into some fixed $2+q$ dimensional background 
that are stationary points of the said functionals under deformations of the 
embedding in any normal direction. 
We may write the mean curvature as 
$H=h\nu_H$, where $h$ is a scalar and $\nu_H$ is a normal vector in the 
direction of $H$ (see details in \cite[p. 16]{scs}), and let  
$A_{\nu}$ be the shape operator along $\nu$, and 
$\nabla^{\nu}$ be the covariant derivative of the normal bundle, respectively.
If $\{ e_1, e_2\}$ and $\{\nu_H,\nu_2, \ldots, \nu_q\}$ are orthonormal frames 
of tangent vectors, and normal bundle, 
a normal stationary critical embedding $f$ of (\ref{fu1}) satisfies the
system of equations    
\begin{equation} \label{cphi}
\begin{array}{rcl}
2\Delta h  & = & h\left( 2(K^{\tg}(e_1,\nu_H) +K^{\tg}(e_2,\nu_H) ) 
-2\| \nabla_{e_i}^{\nu}\nu_H \|^2 -
h^2 +2\, {\rm trace}\, A_{\nu_H}^2\right) \, , \\
0 & = & 2h\< R^{\tg}(\nu_H,e_i)e_i,\nu_m\>+4e_i(h)\< \nabla_{e_i}^{\nu}\nu_H,
\nu_m\>+2he_i \< \nabla_{e_i}^{\nu}\nu_H,\nu_m\>
 \\ & & - 2h\< \nabla_{e_i}^{\nu} \nu_H, \nabla_{e_i}^{\nu} \nu_m \> 
+2h\, {\rm trace}\, A_{\nu_H}A_{\nu_m}\, ,   \; m=2, \ldots, q\, ,    
\end{array} 
\end{equation}
where $R^{\tg}$ is the curvature tensor of $\tg$ \cite[Theorem 3.10]{gracie}. 
On the other hand, the variation of $\Theta_f$ along deformations with 
variational vector $T=T^\tau+T^\nu$ is given by 
$$
\int_{f_g(\Sigma)}\sum K^{\tg}(e_i,\nu_H)({\rm div}\, T^\tau-\< H,T\>)d\mu_g
$$ 
\cite[Theorem 3.12]{gracie}, and so its variation along $T^\nu=\nu_H$ 
is $-h\sum K^{\tg}(e_i,\nu_H)$, and vanishes along $\nu_m$, $m=2,\ldots, q$.
Hence, if $f : (\Sigma,g) \hookrightarrow (\mb{S}^{2+q}, \tg)\subset(
\mb{S}^{\tn},\tg)$ is a critical point of $\mc{W}_{f_g}(\Sigma)$,
 the embedding $f$ must satisfy the system of equations  
\begin{equation}  \label{wcc}
\begin{array}{rcl}
2\Delta h  & = & h\left( -2 \| \nabla^{\nu}_{e_i}\nu_H \|^2 -
h^2 +2\, {\rm trace}\, A_{\nu_H}^2\right) \, , \\
0 & = & 4e_i(h)\< \nabla_{e_i}^{\nu}\nu_H,
\nu_m\>+2h \< \nabla_{e_i,e_i}^2\nu_H,\nu_m\>
+2h\, {\rm trace}\, A_{\nu_H}A_{\nu_m}\, ,   \; m=2, \ldots, \tn\, .
\end{array}
\end{equation}

A minimal $f_g\in \mc{M}_{\mb{S}^{\tn}}(\Sigma)$ is 
a critical point of the Willmore functional,
and a normal stationary critical point of $\Psi_f(\Sigma)$. These 
functionals admit critical points that are nonminimal, but the former admits
deformations of those onto minimal ones, while the latter does not.  
In fact, if the Willmore functional is defined over a restricted background, 
hence a different domain,  
you could then have paths of embeddings $f$s of nonconstant $\| H_f\|^2$ 
satisfying the system (\ref{wcc}) each, 
and which can be continuously deformed within this 
background to a minimal embedding. For if the embedding is 
in codimension $q=1$, (\ref{wcc}) reduces to the single equation
$$
2\Delta h = h(-h^2+2\| \alpha \|^2)= (k_1-k_2)^2 h= (k_1^2-k_2^2)(k_1-k_2) \, ,
$$
$k_1$ and $k_2$ the principal curvatures of the surface.\footnotemark \, 
Pinkal exhibits isometric embeddings of flat metric representatives 
of all conformal 
classes of tori into $(\mb{S}^3,\tg)$ that satisfy this equation \cite{pi}, 
vary smoothly with the conformal class, and  from which we can extract 
continuous paths that begin at the embedding of the Clifford torus, but have 
constant mean curvature functions only at the start. The scale of these tori 
changes, and the Clifford torus is the one of smallest area. None of the 
deformations are normal stationary critical points of $\Psi$ in codimension one.

\footnotetext{This was known to Thomsen \cite{thom} when 
$(\tm,\tg)=(S^3_c,\tg)$, $c=0$, or $1$, cases where the variational formula 
for (\ref{fu1}) is seemingly simpler 
\cite[Corollary 3.4]{gracie} (cf. \cite[Theorem 8]{rss2}), and from 
which the said result for $\mc{W}$ is straightforward, as pointed out already.} 

An isometric embedding can be deformed conformally in 
$\mc{M}_{\mb{S}^{\tn}}(\Sigma)$ to the
embedding on an area preserving Einstein metric in the class. The latter 
can then be changed homothetically and then conformally dilated to another 
Einstein $g$ for which the associated embedding 
$f_g \in \mc{M}_{\mb{S}^{\tn}}(\Sigma)$ is 
minimal, and realizes the invariant $\mc{W}_{f_g}(\Sigma)=
\mc{W}(\Sigma,[g])=4\, {\rm area}_g( \Sigma)$ of the class.  
The operations leading to this $f_g$ 
either preserve or lower the value of the Willmore functional of the starting 
metric. If this $f_g$ is a global minimizer of $[g]\rightarrow 
\mc{W}(\Sigma,[g])$, then $\mc{W}_{f_g}(\Sigma)$ is stationary along all 
isometric embedding deformations associated to deformations of $g$ in 
directions orthogonal
to $g \in [g]$. When this is not an empty condition, we conclude
that the fundamental domain of this $g$ is equilateral of
equal interior angles, and the absolute minimum of $\mc{W}$ is given by
$4\, {\rm area}_g(\Sigma)$.
This summarizes our approach to the GWC.

Our main result is the following theorem. We shall prove 
it in \S \ref{s3}, after the discussion of technical preliminaries 
in \S \ref{s2} required for the task.

\begin{theorem} \label{th1}
If $\Sigma=\Sigma^k$ is oriented, and $\mc{W}$ is the Willmore
functional {\rm (\ref{will})}, then:
\begin{enumerate}[label={\rm \alph*)}] 
\item If $k=0$, for any Riemannian metric $g'$ on $\Sigma$, we have that  
$$
\mc{W}_{f_{g'}}(\Sigma)\geq \mc{W}_{f_g}(\mb{S}^2) =16\pi\, ,
$$
where $f_g : (\mb{S}^2,g) \rightarrow (\mb{S}^3,\tg)\subset (\mb{S}^{\tn},\tg)$
is the standard totally geodesic embedding of the two sphere, and $g'\in [g]$.  
The lower bound is achieved by $f_{g'}$ if, and only if, $f_{g'}$ is
conformally equivalent to $f_g$ in $\mb{S}^{\tn}$ and ${\rm area}_{g'}(\Sigma)
\leq 4\pi$.  
\item If $k>0$, for any Riemannian metric $g'$ on $\Sigma$, we have that  
$$
\mc{W}_{f_{g'}}(\Sigma) \geq 
\mc{W}_{f_{g_{\xi_{k,1}}}}(\xi_{k,1}) = 4\, \mu_{g_{\xi_{k,1}}}(\xi_{k,1}) \, ,
$$
and the lower bound is achieved by 
$f_{g'} : (\Sigma , g') \rightarrow (\mb{S}^{\tn},\tg)$ if, and only if,  
$f_{g'}$ 
is conformally equivalent in $\mb{S}^{\tn}$ to the isometric embedding 
$f_{g_{\xi_{k,1}}}: (\xi_{k,1},g_{\xi_{k,1}}) \rightarrow (\mb{S}^3,\tg)
\hookrightarrow (\mb{S}^{\tn},\tg)$ 
and $\mu_{g'}(\Sigma)\leq \mu_{g_{\xi_{k,1}}}(\xi_{k,1})$.
\end{enumerate}
\end{theorem}

Our argument leads to a new proof of the following classic result.

\begin{corollary} \label{co2}
If $\Sigma$ is a Riemann surface of genus zero, then $\Sigma$ is 
$\mb{P}^1(\mb{C})=(\mb{S}^2,\tg)$ as a complex manifold.
\end{corollary}

We obtain also a new proof of the following conjecture of Lawson \cite{law2}, 
a recent theorem of Brendle \cite{bren}. 

\begin{corollary} \label{co3}
If $\Sigma$ is a Riemann surface of genus one minimally embedded into 
$\mb{S}^3$, then up to isometries of the ambient space, $\Sigma$ is the 
Clifford torus.
\end{corollary}

We prove the analogous version of this theorem for nonorientable surfaces in 
\S \ref{s4}.

\section{Conformal deformations of isometric embeddedings of surfaces}
\label{s2}
We denote by $\mc{M}_a(\Sigma)=\{ g\in \mc{M}(\Sigma): \; {\rm area}_g(
\Sigma)=a \}$ the manifold of Riemannian metrics on $\Sigma$ of area $a$,
by $\mc{M}_{a,[g]}(\Sigma)$ the submanifold of metrics of area $a$ in the 
conformal class $[g]$, and by $\mc{M}_{[g]}(\Sigma)$ the submanifold of metrics
in $[g]$. They correspond to the spaces 
$\mc{M}_{\mb{S}^{\tn},a}(\Sigma)=\{ f_g \in \mc{M}_{\mb{S}^{\tn}}(\Sigma)
\, , \; {\rm area}_g(\Sigma)=a\}$, 
$\mc{M}_{\mb{S}^{\tn},a,[g]}(\Sigma)=\{ f_g\in \mc{M}_{\mb{S}^{\tn},a}(\Sigma):
g\in [g]\}$, and 
$\mc{M}_{\mb{S}^{\tn},[g]}(\Sigma)=\{ f_g\in \mc{M}_{\mb{S}^{\tn}}(\Sigma):
g\in [g]\}$, respectively. If convenient, we shall write ${\rm area}_g(\Sigma)$
as $\mu_g(\Sigma)$, and go back and forth between the two notations.
\smallskip

We recall that any metric $g\in \mc{M}_{a,[g]}(\Sigma)$ can be conformally
deformed in $\mc{M}_{a,[g]}(\Sigma)$ to an extremal $g^e$ of the 
functional  
\begin{equation} \label{ss}
\mc{S}_{a,[g]}^2: \mc{M}_{a,[g]}(\Sigma) \ni g \rightarrow \int s_g^2 d\mu_g\, ,
\end{equation}
an Einstein minimizer of scalar curvature $4\pi \chi(\Sigma)/a$, 
$\chi(\Sigma)$ the Euler characteristic of $\Sigma$. 
This metric $g^e$ is unique modulo isometries, and varies differentiably 
with $[g]$. If $\Sigma$ is oriented, the extremals of $\mc{S}_{a,[g]}^2$ 
are just the extremal metrics of Calabi \cite{ca} on the K\"ahler curve 
$(M,J,g)$, $J$ the complex structure induced by $g$. They can be solved for by 
deforming the metric $g$ along a path  
$t \mapsto  \omega_{g_t}=\omega_g + i \ddb \varphi(t)$ 
of K\"ahler forms in the class $\Omega_g= [\omega_g]$, 
$\omega_g(\, \cdot \, ,\, \cdot \,)= g(J\, \cdot \, ,\, \cdot \,)$,
evolving in time according to the Yamabe flow of Hamilton \cite{haY}. 
(For dimensional reasons, this and the extremal 
flow \cite{simH} coincide.) On a Riemann surface, 
this flow equation has a solution for all time that converges at infinity 
to a metric $g^e$ of constant scalar curvature $s_g$ minimizing 
the functional, the value of $s_g$ following by the extreme 
case of the Cauchy-Schwarz inequality, and the Gauss-Bonnet 
theorem\footnotemark. \hspace{-2mm}
By Moser's theorem \cite{mos}, there exists  
an associated family of symplectomorphims $\eta(t): (\Sigma,\omega_g) 
\rightarrow (\Sigma,\omega_{g_t})$ realizing the conformal deformations of 
the metrics, $\omega_{g_t} = e^{2\psi(t)}\omega_g$ for a 
suitable conformal factor $e^{2\psi(t)}$, and $\eta(t)^* \omega_{g_t}=
\omega_g$. Since the function $[g] \rightarrow J$ is
regular, this area $a$ metric $g^e$ varies differentiable with the
class. If $\Sigma$ is nonoriented, this argument can be carried out on the 
oriented 
double cover with the lifted metric. The push down of the solutions to 
the flow equation there produces a path of area preserving 
metrics on the base, all conformally equivalent to $g$, with the desired
properties, the areas of the metrics in cover and base differing by a factor 
of $1/2$. The diffeomorphism realizing the conformal relation of 
the metrics upstairs does not descend to a likewise diffeomorphism 
on the base, the 
area forms of the conformally related metrics on the base being a 
well-defined $1$-density only.
\smallskip
\footnotetext{This provides an alternate proof of the uniformization theorem 
for Riemann surfaces.  For $\mb{S}^2$, Hamilton's result was short of the 
general statement, but was completed later with a key contribution by B. Chow 
\cite{cho}, and an additional observation by X.X. Chen, P. Lu and G. Tian 
\cite{xxluti}.} 

We study next the isometric embeddings into a fixed background of conformally 
related metrics on $M$. Given the embedding (\ref{emb}) of $(M,g)$ into 
$(\tm,\tg)$, suppose that 
\begin{equation} \label{emc}
[0,1]\ni t \rightarrow f_{g_t}: (M,g_t) \mapsto (\tm, \tg)
\end{equation}
is a path of conformally related isometric embeddings deformations of 
$f_g=f_{g_t\mid_{t=0}}$ of conformally related Riemannian metrics 
\begin{equation} \label{cd}
[0,1]\ni t \rightarrow g_t =e^{2\psi (t)}g
\end{equation}
on $M$, $\psi(t)$ a path of functions, $\psi(0)=0$.
Under mild conditions on $(\tm,\tg)$, the Palais isotopic extension 
theorem \cite{pal} applies, and we obtain a smooth one parameter family of 
diffeomorphisms
$$
F_t: \tm \rightarrow \tm 
$$
such that
$F_t(f_g(x))=f_{g_t}(x)$. 
Since the metrics on the submanifolds are all induced by
the metric $\tg$ on the background $\tm$, we obtain by restriction a 
diffeomorphism 
$$
F_t: (f_g(M),\tg) \rightarrow (f_{g_t}(M), \tg) \, . 
$$
The tensor $F_t^* \tg$ is just the metric $\tg$ on $\tm$ acted on by the
diffeomorphism $F_t$, and since the $g_t$s are conformal deformations of
$g$, we have that
$$
F_t^* (\tg\mid_{f_{g_t}(M)}) = e^{2 u(t)(f_g(\, \cdot \,))} \tg \mid_{f_g(M)} =
e^{2 u(t)(f_g(\, \cdot \,))} g  \, ,
$$
where the conformal factor $e^{2 u(t)}$ and that in (\ref{cd}) are related 
to each other by $e^{2\psi(t)(\, \cdot \, )}=e^{2u(t)\circ f_g(\, \cdot \, )}$. 
We extend $u(t)$ conveniently to a function on the whole of
$\tm$, and view the family (\ref{emc}) as the family of conformal isometric 
embeddings 
\begin{equation}\label{embt}
f_{g_t}: (M,e^{2u(t)\circ f_g}g) \hookrightarrow (\tm ,e^{2u(t)}\tg)
\end{equation}
of the fixed manifold $f_g(M)$ with a varying metric on it. We 
compute accordingly the densities of the functionals 
$\Theta$, $\Psi$ and $\Pi$ associated to $f_{g_t}$ in terms of the densities 
of these functionals associated to $f_g$. 
For notational simplicity, we set $\tilde{g}_t=e^{2u(t)}\tilde{g}$.

If $e_1, \ldots, e_n$ is a $\tilde{g}$-orthonormal tangent basis on 
$f_g(M)$, 
then $e_1^t=e^{-u(t)}e_1, \ldots, e_n^t=e^{-u(t)}e_n$ is an orthonormal 
tangent basis for $f_{g_t}(M)$, and so its exterior scalar curvature transforms
as 
\begin{equation} \label{eq1}
\sum_{i,j} K^{\tilde{g}_t}(e_i^t, e_j^t)= e^{-2u(t)}(
\sum_{i,j} K^{\tilde{g}}(e_i,e_j) -2(n-1)(
{\rm div}_{f_g(M)} \nabla^{\tg} u^\tau- \tg(H, \nabla^{\tg}
u^\nu)-|du^\tau|^2+\frac{n}{2}|du|^2 )) \, , 
\end{equation}
where $\tau$ and $\nu$ denote the tangential and normal components, 
respectively. On the other hand, by the relationship between 
the second fundamental forms of $f_{g_t}(M)$ and $f_g(M)$, 
we have that
\begin{eqnarray}
\| H \|^2_{f_{g_t}} & = & e^{-2u(t)}(\| H\|^2_{f_g} 
- 2n\tg(H,\nabla^{\tg}u ^\nu)
+n^2 \tg(\nabla^{\tg}u ^\nu, \nabla^{\tg} u^\nu)) \, , \label{eq3}\\
\| \alpha\|^2_{f_{g_t}} & = & e^{-2u(t)}(\| \alpha\|^2_{f_g} - 2\tg(H,
\nabla^{\tg}u^\nu) +n \tg(\nabla^{\tg}u ^\nu, \nabla^{\tg} u^\nu)) \, .
\label{eq4} 
\end{eqnarray}
These three expressions on $f_g(M)$ are fully determined once we 
know the first jet in the normal directions of the extension of $u(t)$ to a 
tubular neighborhood of $f_g(M)\subset \tm$, and the last two are independent
of the first one. 
By (\ref{sce}), they imply the intrinsic scalar curvature relation 
\begin{equation} \label{eq5}
s_{g_t}= e^{-2u(t)}\left( s_{g}-2(n-1){\rm div}_{f_g(M)}(\nabla^{\tg}u)^\tau
-(n-1)(n-2)\tg(\nabla^{\tg}u^{\tau},\nabla^{\tg}u^\tau)  \right) \, .
\end{equation}

If $\Sigma=M^2$ is two dimensional, and each of (\ref{eq1}), (\ref{eq3})
and (\ref{eq4}) is satisfied, 
\begin{eqnarray}
\Theta_{f_{g_t}}(\Sigma) & = & \Theta_{f_g}(\Sigma)+2
\int_{f_g(\Sigma)} (\< H , \nabla u^\nu\> - \< \nabla u^\nu, \nabla u^\nu\>)
 d\mu_g \, , \label{eqn1} \\
\Psi_{f_{g_t}}(\Sigma) & = & \Psi_{f_g}(\Sigma)-4\int_{f_g(\Sigma)}(\< H, 
\nabla u^\nu\> - \< \nabla u^\nu, \nabla u^\nu\>) d\mu_g \, , \label{eqn2} \\ 
\Pi_{f_{g_t}}(\Sigma) & = & \Pi_{f_g}(\Sigma)-2 \int_{f_g(\Sigma)} (\< H, 
\nabla u^\nu\> - \< \nabla u^\nu, \nabla u^\nu\>) d\mu_g  \, ,
\label{eqn3}
\end{eqnarray}
respectively, and we have that
$$ 
\mc{W}_{f_{g_t}}(\Sigma) = 2\Theta_{f_{g_t}}(\Sigma)+\Psi_{f_{g_t}}(\Sigma) = 
2\Theta_{f_g}(\Sigma)+\Psi_{f_g}(\Sigma) =\mc{W}_{f_g}(\Sigma)  \, ,   
$$ 
indicating that the Willmore value of the starting $f_g$ constraints all the 
values of $\mc{W}_{f_{g_t}}$ when the $f_{g_t}$s are conformally related to
each other. Notice that if
$$
\mc{D}_{f_g}(\Sigma)= \Theta_{f_g}(\Sigma) + \Pi_{f_g}(\Sigma)\, ,
$$
and we proceed similarly, we have that
$$
\mc{D}_{f_{g_t}}(\Sigma) = \Theta_{f_{g_t}}(\Sigma)+\Pi_{f_{g_t}}(\Sigma) =
\Theta_{f_g}(\Sigma)+ \Pi_{f_g}(\Sigma)= \mc{D}_{f_g}(\Sigma)\, ,
$$
and the total scalar curvature (\ref{tsc}) of $(\Sigma, g)$   
becomes the difference
$$
\mc{S}_g(\Sigma)= 4\pi \chi(\Sigma)= \mc{W}_{f_g}(\Sigma)-\mc{D}_{f_g}(\Sigma) 
$$
of functionals invariant under conformal deformations $f_{g_t}$ of
$f_g$. 

\begin{theorem} \label{th4}
If $f_{g_t}$ is a path of conformal diffeomorphism deformations of an 
isometric embedding $f_g: (\Sigma,g)\rightarrow (\tm,\tg)$ of a 
surface $(\Sigma,g)$, then $\mc{W}_{f_{g_t}}(\Sigma)=\mc{W}_{f_g}(\Sigma)$, 
and along such deformations, if 
any one of the three functionals $\Theta_{f_g}(\Sigma)$, $\Psi_{f_g}(\Sigma)$, 
or $\Pi_{f_g}(\Sigma)$ is preserved, so are the other two.
In particular, if $(\tm,\tg)=(\mb{S}^{\tn},\tg)$ and 
$F: (\mb{S}^{\tn},\tg) \rightarrow (\mb{S}^{\tn},\tg)$ is a 
conformal automorphism, 
$\mc{W}_{f_g}(\Sigma)=\mc{W}_{F\circ f_g}(\Sigma)$. 
\end{theorem}

{\it Proof}. If along the deformations $f_{g_t}$ any one of 
$\Theta_{f_g}(M)$, $\Psi_{f_g}(M)$, or $\Pi_{f_g}(M)$, 
is preserved, then by either (\ref{eqn1}), (\ref{eqn2}), or (\ref{eqn3}), we 
conclude that 
$$ 
\int_{f_g(\Sigma)} (\< H, \nabla u^\nu\> -
\< \nabla u^\nu, \nabla u^\nu\> )d\mu_g =0 \, , 
$$ 
and so $\Theta_{f_{g_t}}(\Sigma)= \Theta_{f_g}(\Sigma)$, 
$\Psi_{f_{g_t}}(\Sigma)= \Psi_{f_g}(\Sigma)$, and $\Pi_{f_{g_t}}(\Sigma) = 
\Pi_{f_g}(\Sigma)$, respectively. 

The last statement follows by  realizing $F$ as the end point of a path
of conformal diffeomorphisms of the background sphere 
starting at $\BOne_{\mb{S}^{\tn}}$. 
\qed 

\begin{theorem} \label{th5}
The restriction of the Willmore functional to 
$\mc{M}_{\mb{S}^{\tn},a,[g]}(\Sigma)$ is the constant 
$\mc{W}_{f_{g^e}}(\Sigma)$, $f_{g^e}$ the Nash isometric 
embedding of the unique {\rm (}up to isometries{\rm )} Einstein $g^e$ in 
$\mc{M}_{a,[g]}(\Sigma)$. If $f_{g_t}$ is a path in
$\mc{M}_{\mb{S}^{\tn},a}(\Sigma)$, and $f_{g_t^e} \in
\mc{M}_{\mb{S}^{\tn}a,[g_t]}(\Sigma)$ is the corresponding path of isometric
embeddings of the Einstein $g_t^e$s in the $[g_t]$s, we have that
$\mc{W}_{f_{g_t}}(\Sigma)=\mc{W}_{f_{g^e_t}}(\Sigma)$.
\end{theorem}

{\it Proof}. We consider $g \in \mc{M}_{a,[g]}(\Sigma)$ with associated
Nash isometric embedding $f_g$, and let $f_{g_t}$ be a path in 
$\mc{M}_{\mb{S}^{\tn},a,[g]}(\Sigma)$ of deformations of $f_g$ corresponding 
to a path $g_t \in \mc{M}_{a,[g]}(\Sigma)$ of conformal deformations of 
$g$. If $u=u(t)$ is the path of scalar functions in (\ref{embt}) realizing 
the conformal relation $\tilde{g}_t =e^{2u(t)} \tilde{g}$, by (\ref{eq1}), 
$u$ satisfies the strong equation 
$$
2 = e^{-2u} (2+2 \Delta u +2(\< H, \nabla^{\tg}u^\nu\>- \< \nabla^{\tg}u^\nu,
\nabla^{\tg}u^\nu\>))
$$
for each $t$. By integration relative to the area measures of the area
preserving conformally deformed metrics, we conclude that 
$$ 
\int_{f_g(\Sigma)} (\< H, \nabla u^\nu\> -
\< \nabla u^\nu, \nabla u^\nu\> )d\mu_g =0 \, ,  
$$ 
so $\Theta_{f_g}(\Sigma)$ is preserved along $f_{g_t}$,    
and by Theorem \ref{th4}, $\mc{W}_{f_{g_t}}(\Sigma)=\mc{W}_{f_g}(\Sigma)$.
In particular, if we consider the path $g_t$ to the Einstein minimizer 
$g^e$ of the functional (\ref{ss}), 
we conclude that $\mc{W}_{f_{g^e}}(\Sigma)=\mc{W}_{f_g}(\Sigma)$.

Given $f_{g_t}\in \mc{M}_{\mb{S}^{\tn},a}(\Sigma)$, we take the path of 
intrinsic metrics $g_t\in \mc{M}_a(\Sigma)$ as a smooth family of initial 
conditions for the Yamabe flows on the conformal classes $[g_t]$s. The 
solutions produce a path of extremal minimizers $g^e_t$ of (\ref{ss}) in 
each $[g_t]$, with $s_{g^e_t}=4\pi \chi(\Sigma)/a$ each. By the previous 
result, $\mc{W}_{f_{g_t}}(\Sigma)=\mc{W}_{f_{g_t^e}}(\Sigma)$. 
\qed

The result above elucidates the behaviour of $\mc{W}$ on isometric embeddings 
of area preserving metrics in a given conformal class. By using 
the variational expression
\begin{equation} \label{af}
\frac{d}{dt}d\mu_{g_t}=\left( {\rm div}(T^{\tau}) - \< T^{\nu},
H_{f_{g_t}}\>\right) d\mu_{g_t}\, .
\end{equation}
for the volume forms along $f_{g_t}$,
$T=T^{\tau}+T^{\nu}$ the variational vector field, combined with (\ref{eq3}),
we can do likewise for the ray of homothetics of these metrics.

\begin{theorem} \label{th6}
The quantity {\rm (}\ref{min}{\rm )} is a well defined conformal invariant of 
$[g]$, and if $g^e \in [g]$ is Einstein, it can be computed by
$$
\mc{W}(\Sigma,[g])= \inf_{s\in \mb{R}}\mc{W}_{f_{e^{2s}g^e}}(\Sigma)\, .
$$
The function $[g] \rightarrow \mc{W}(\Sigma,[g])$ is differentiable.
\end{theorem}

{\it Proof}. Let $g\in [g]$, and consider the conformally related metric
$g'=e^{2u}g$. By the transformation rule of the area forms, and (\ref{eq3}), 
we have that 
$$
\begin{array}{rcl}
\mc{W}_{f_{e^{2v}g'}}(\Sigma) & = &  {\displaystyle 
\int_{f_{e^{2v}g'}(\Sigma)} ( 4 + \| H_{f_{e^{2v}g'}}\|^2 )
d\mu_{f_{e^{2v}g'}} } \\ & = & {\displaystyle   
\int_{f_{g}(\Sigma)} e^{2(u+v)}( 4 + \| H_{f_g}-2\nabla (u+v)^{\nu}\|^2 ) 
d\mu_{f_{g}} }  \\ & = &
\mc{W}_{f_{e^{2(u+v)}g}}(\Sigma)\, ,
\end{array}
$$
and the infima of the starting and ending expressions coincide.

By Theorem \ref{th5}, the invariant $\mc{W}(\Sigma,[g])$ may be computed
by taking the infimum over the isometric embeddings of the homothetics of 
$g^e$, the area preserving Einstein deformation of $g$ in the class $[g]$. 
If we let $\Sigma_j=f_{e^{2t_j}g^e}(\Sigma)$ be a minimizing sequence,
we may use translations and scalings to assume that $p=(1,0, \ldots, 0)\in 
\Sigma_j$, and that $\mu_{e^{2t_j}g^e}(\Sigma)\geq 1$. Thus, the sequence of 
Gaussian curvatures is bounded, and by \cite[Lemma 1.1]{lsim}, the sequence of 
diameters is bounded above. This places us in the context Gromov-Hausdorff 
convergence of pointed spaces with suitable bounds. The argument of 
\cite[Theorem 3.1]{lsim} may be applied in this restricted case, and we 
obtain a subsequence $\Sigma_{j'}$ and a real analytic surface $\mc{S} \subset
\mb{S}^{\tn}\subset \mb{R}^{\tn +1}$ such that $\Sigma_{j'}$ converges to
$\mc{S}$ in the Gromov-Hausdorff distance, and achieves $\mc{W}(\Sigma,[g])$. 

The structure of the limit $\mc{S}$ may be argued in detail when $\Sigma$ is
oriented. In that case, the same idea as in \cite[\S 4]{lsim} applies for the 
cases of genus  $k=0,1$, to conclude that $\mc{S}$ must a sphere or a torus,
respectively. If $k\geq 2$, by Bishop inequality and the Gauss-Bonnet theorem, 
we can prove that diameters of $\Sigma_j$ cannot go to zero, and so 
$\mc{S}$ must have dimension two also, and it is an oriented surface of 
genus $k$. In general, we may take advantage of the conformal relation
of the Einstein metrics to draw the appropriate conclusion, since we already
know a limit exists. If $f_{g^e}$ is not minimal to begin with, 
by (\ref{eq3}), if we take a sequence of deformations in the 
conformal direction $H_{f_{g^e}}$, the value of $\| H\|^2_{f_{g^e}}$ decreases.
If these deformations are conformal dilations so by (\ref{eq1}) the area 
increases, by Theorem \ref{th4}, the Willmore energy along them is preserved; 
otherwise, by (\ref{af}), the area of the embedded deformations decreases also, 
and the Willmore energy decreases. Hence, the infimum is achieved at a 
homothetic deformation $f_{e^{2s}g^e}$ that is minimal, or at any 
$f_{e^{2s}g^e}$ for which deformations of the metric along 
$H_{f_{e^{2s}g^e}}$ are conformal dilations of increasing area, 
in either case, the intrinsic metric of the limit being Einstein. 

If $g_t$ is a path of area preserving deformations of $g$, and $g_t^e$
is the corresponding path of area preserving Einstein deformations of 
it in each class $[g_t]$, we have that
$$ 
\mc{W}(\Sigma,[g_t])= \inf_{s_t \in \mb{R}}\mc{W}_{f_{e^{2s_t}g_t^e}}(\Sigma)
$$ 
where $e^{2s_t}g_t^e$ is the homothetic of $g_t^e$, and $f_{e^{2s_t}g_t^e}$ 
is its associated isometric embedding, respectively. This function is
differentiable at $[g_0]=[g]$, hence everywhere.
\qed
\smallskip

\begin{theorem} \label{th7}
Suppose that 
$f_g: (\Sigma,g) \rightarrow (\mb{S}^{2+q},\tg) \hookrightarrow (\mb{S}^{\tn},
\tg)$ is a minimal embedding of codimension $q$. If $f_{g_t}\in 
\mc{M}_{\mb{S}^{\tn}}(M)$ is a path of embedding deformations of $f_g$ 
associated to a path $g_t$ of area preserving conformal deformations of $g$, 
then each of the embeddings $f_{g_t}$ is minimal, of codimension $q$ also,
and $\mc{W}_{f_{g_t}}(\Sigma)=4\, {\rm area}_g(\Sigma)$. If $f_g$ and
$f_{g'}$ are minimal isometric embeddings of metrics $g,g'$ in the same 
conformal class, then ${\rm area}_{g}(\Sigma)={\rm area}_{g'}(\Sigma)$.
\end{theorem}

{\it Proof}. We let $u=u(t)$ be the path of scalar functions above such that 
$\tilde{g}_t =e^{2u(t)} \tilde{g}$. By (\ref{eq1}), for each $t$, 
$u$ satisfies the strong equation   
$$ 
2=e^{-2u(t)}(2+2\Delta^{g}u 
-2\< \nabla^{g} u^\nu, \nabla^{g} u^\nu )\> \, . 
$$ 
By integration relative to the area measures of the conformally deformed 
metrics,
we conclude that the global $L^2$-norm of $\nabla u^{\nu}$ is zero, and so 
$u(t)$ does not vary in the normal directions. Thus, the codimension of all the 
embedding deformations remains constant, and by (\ref{eq3}), these 
embedding deformations are all minimal. 

If $f_g$ and $f_{g'}$ are minimal, by Theorem \ref{th5} and this preliminary
result, we can assume that $g$ and $g'$ are Einstein, and so
a homothetic of each other.  If $[0,1]\ni t \rightarrow \tilde{g}_t=
e^{2u}\tilde{g}$ is a path of conformal deformations of $g$ into $g'$ as above, 
by (\ref{eq3}), 
$$
\| H\|^2_{f_{g'}}= 4 e^{-2u} \< \partial u^{\nu}, \partial u^{\nu}\>\mid_{t=1}
\, ,
$$
and so $\< \partial u^{\nu}, \partial u^{\nu}\>\mid_{t=1}=0$. By (\ref{eq1}),
$g'=e^{2u\mid_{t=1}}g$ and $g$ must be metrics of the same area.
\qed
\smallskip

By Theorem \ref{th4}, a minimal critical point $f_g$ of the Willmore
functional $\mc{W}_{f_g}(\Sigma)$ admits conformally related continuous 
deformations onto nonminimal ones that maintain the critical value.
This contrasts with the normally stationary critical points of 
$\Psi_{f_g}(\Sigma)$, and becomes clear when we 
cast the gap theorem of Simons \cite{si,cdck,law0} as a statement for these 
type of critical embeddings of $\Psi_f(M)$ \cite[Theorem 2]{rss2}.
For suppose that $f_g \in \mc{M}_{\mb{S}^{\tn}}(\Sigma)$ is a 
critical point of $\Psi_{f}(\Sigma)$ of codimension $q$, so $f_g$ satisfies 
the system (\ref{cphi}). If the first of the 
inequalities in 
\begin{equation} \label{es}
\begin{array}{rcl}
{\displaystyle -\lambda \| H\|^2 -2} & \leq &
{\displaystyle {\rm trace}\, A_{\nu_H}^2-\| H\|^2 -\| \nabla^\nu \nu_H \|^2} \vspace{1mm} \\
 & \leq & {\displaystyle \| \a \|^2 -\| H\|^2 -\| \nabla^\nu \nu_H \|^2}
\leq {\displaystyle \frac{2q}{2q-1}}
\end{array}
\end{equation}
holds for some $\lambda \in [0,1/2)$, then $f_g$ is minimal. In that case,
if 
\begin{equation} \label{esu1}
0 \leq \| \alpha \|^2 \leq 2q/(2q-1) \, ,
\end{equation}
then either $\| \alpha \|^2 = 0$ and $(\Sigma,g) = (\mb{S}^2,\tg)$, or 
$\| \alpha \|^2 = 2q/(2q-1)$, and either $\| \alpha \|^2=2$, and $(\Sigma,g)$ 
is $(\xi_{1,1}, g_{\xi_{1,1}})$, or $\| \alpha \|^2=4/3$ and $(\Sigma,g)$ is 
$\mb{P}^2(\mb{R})$ minimally embedded into $\mb{S}^4$ with an Einstein metric 
of scalar curvature $2/3$. Otherwise,  
\begin{equation} \label{esu2}
\| \alpha \|_{f_g}^2 > 2q/(2q-1)
\end{equation}
at at least one point of $f_g(\Sigma)\subset \mb{S}^{2+q}\subset \mb{S}^{\tn}$.
If, on the other hand, the first of 
the inequalities in (\ref{es}) fails to hold for any $\lambda \in
[0,1/2)$, then $f_g$ cannot be minimal, and the first of the equations in
(\ref{cphi}) implies that 
$$
0 \leq {\rm trace}\, A_{\nu_H}^2=
\frac{1}{2}\| H\|_{f_g}^2 + \| \nabla^{\nu}\nu_H\|^2 -2 \leq 
\| \alpha \|^2_{f_g} \, ,  
$$ 
which, by the nonnegativity of the left side, shows that this type of 
critical embedding cannot be continuously deformed onto a minimal one.
The embeddings of the round two sphere in $\mb{S}^3$ as a sphere of radius $r$
illustrates this dichotomy well. 

The dichotomy above remains so even if we allow  
for continuous variations of the conformal class $[g]$ of the critical point,  
and corresponding variations of the critical values $\mc{W}(\Sigma,[g])$. In
these cases, the only paths that are also critical points of $\Psi$ are
paths of minimal embeddings across the different classes. Even the flat tori by 
Pinkall \cite{pi} alluded to earlier exhibits this property, though they
are critical points of $\mc{W}_{f_g}$ in their conformal classes of 
nonconstant $\| H\|^2$.  
\medskip

\begin{theorem} \label{th8}
If $f_g\in \mc{M}_{\mb{S}^{\tn}}(\Sigma)$ is such that $\mc{W}_{f_g}(\Sigma)=
\mc{W}(\Sigma,[g])$, and $a={\rm area}_g(\Sigma)$, then $f_g$ is a critical 
point of the functional
\begin{equation}
\label{ssss}
\mc{W}_a: \mc{M}_{\mb{S}^{\tn}, a}(\Sigma) \ni f_{g'}
\rightarrow \mc{W}_{f_{g'}}(\Sigma)
\end{equation}
if, and only if, the function $g' \rightarrow {\rm area}_{g'}(\Sigma)$ 
is stationary along deformations of $g\in [g]$ that are 
orthogonal to $g$. Thus, if $\Sigma^k$ is oriented of genus $k \geq 1$, 
the fundamental domain of $(\Sigma, [g])$ is an equilateral 
geodesic polygon of interior angles all equal to $2\pi/4k$, 
and the critical value is $\mc{W}_{f_g}(\Sigma)=(4+\| H_{f_g}\|^2)a$,
while if $\Sigma^k$ is nonoriented of genus $k\geq 2$, 
the lifted metric $\overline{g}$ of $g$ to the 2-to-1 oriented cover 
$\tilde{\Sigma}^{k-1}$ must be in the conformal class of metrics with 
this type of fundamental domain, in which case the
critical value is $\mc{W}_{f_g}(\Sigma)=(4+\| H_{f_{g}}\|^2)a$ also.
\end{theorem}

{\it Proof}. Since $f_g$ realizes the invariant $\mc{W}(\Sigma,[g])$ of the
conformal class, the variations of $\mc{W}_{f_g}(\Sigma)$ in 
$\mc{M}_{\mb{S}^{\tn},[g]}(\Sigma)$ vanish. This reduces
our considerations to variations of $\mc{W}_a$ along deformations 
$f_{g_t}$ of $f_g$ corresponding to variations $g_t$ of 
$g$ in directions orthogonal to $g \in [g]$ itself, that is to say, along
symmetric two tensors of zero $g$ trace. By (\ref{af}),  
the initial deformation vector $T$ of these deformations is tangent to 
$f_g(\Sigma)$. We compute the variations of $\mc{W}_a$ along this type of 
deformations.

We proceed as in \cite{gracie}. We take an 
arbitrary point $p\in f_g(\Sigma)$, and on a sufficiently small neighborhood 
of the integral curves of the deformation vector $T$ through 
$(p,0) \in f_{g}(\Sigma) \times \{ 0\}$, we choose an orthonormal tangent 
frame $\{e_1,e_2\}$ that is geodesic, and commutes with $T=T^\tau$ at the 
said point. Since 
$$
\| H_{f_{g_t}}\|^2 =
\| h_{f_{g_t}}\nu_{H_{f_{g_t}}}\|^2 =
 g^{ij}g^{lk}\< \a_{f_{g_t}}(e_i,e_j),\a_{f_{g_t}}(e_l,e_k)\>\, ,
$$  
dropping the metric from the notation when $t=0$, we have that 
$$
\begin{array}{rcl}
{\displaystyle \frac{d}{dt}\< H_{f_{g_t}},H_{f_{g_t}}\>\mid_{t=0}} & = &   
2h^2\< \nabla_T^{\tg} \nu_H,\nu_H\>+2\dot{g}^{ij}\< \a (e_i, e_j), H)\>  
\\ & = & 
4 \< \a (e_i,e_j),H\>\< \a (e_i,e_j),T\> -
2 (\< e_i\< T, e_j\>+e_j \< T, e_i\>)\a (e_i,e_j),H\> \vspace{1mm}\\ & = & 
-4 \< \a (e_i,\nabla_{e_i}^g T),H\> \\  & = & 0 \, ,   
\end{array}
$$
cf. \cite[Theorem 3.2]{gracie}.
Thus, the function ${\displaystyle \frac{d}{dt}\|H_{f_{g_t}}\|^2\mid_{t=0}}$ 
vanishes identically, and we conclude that
$$
{\displaystyle \frac{d}{dt}\mc{W}_{f_{g_t}}(\Sigma)\mid_{t=0}} =  
{\displaystyle \frac{d}{dt}(4 \,{\rm area}_{g_t}(\Sigma) + \| H_{f_{g_t}}\|^2 
{\rm area}_{g_t}(\Sigma))\mid_{t=0}} 
={\displaystyle (4+\|H_{f_g}\|^2)\frac{d}{dt}\,{\rm area}_{g_t}(
\Sigma)\mid_{t=0}}\, ,
$$
proving the stated characterization of the critical point $f_g$ of 
$\mc{W}_a$.

By Theorem \ref{th5}, we may assume that the 
metric $g$ at the start is Einstein. We proceed with the proof assuming 
firstly that $\Sigma$ is oriented. We let $D$ be the fundamental domain $D$ of 
the class $[g]$ that $g$ represents. We extremize ${\rm area}_g(\Sigma)$ in 
terms of its conformal parameters, which up to a scaling factor, are the
extremals of the normalized ${\rm area}(D)$. 

If $k=1$, for any conformal class of metrics on $\Sigma$ there exists a 
unique parameter $\tau= l e^{i\theta}$, $l|\cos{\theta}|\leq \frac{1}{2}$, 
$l\geq 1$, and the class is represented by the torus $T^2_{1,\tau}$ quotient 
of $\{ z\in \mb{C}: {\rm Im}\, z >0\}$ by the lattice $\<1,{\tau}\>$ of 
rank two with generators $\{ 1, \tau \}$. The area of $D$ for this torus is 
$$
{\rm area}(D)= l\sin{\theta} \, , 
$$ 
which we assume normalized to be $1$. Then the extremal occurs when the 
$\theta=\pi/2$, and $l=1$.

If $k>1$, the fundamental domain $D$ is a geodesic $4k$gon whose vertices
all lie in an orbit of the Fuchsian group of the surface. If $\{ \alpha_i\}$ 
is the set of interior angles of $D$, we have that
$$
{\rm area}(D) = 4k \pi - \sum \alpha_i - 2\pi \, ,
$$
with the angles subjected to the constraints $\alpha_i \geq 0$ and 
$\sum \alpha_i = 2\pi$. By the Gauss-Bonnet theorem, this expression is 
normalized and equal $2\pi |\chi(\Sigma)|$, and the desired extremal of the
linear problem with linear boundary constrains happens at the corner given 
by $\alpha_i = 2\pi/4k$ for all $i$. 

The assertion for the critical value $\mc{W}_{f_g}(\Sigma)$ of the critical 
$f_g$ follows by definition. 

If $\Sigma^k$ is not oriented, its orientable $2$-to-$1$ cover 
$\tilde{\Sigma}$ is a surface of genus $k-1$ with a smooth fixed point free 
involution that defines the cover, and any metric on $\Sigma$ can be lifted 
to a metric on the cover that has a $\mb{Z}/2$ invariance. If $f_g$ is a 
critical point of $W_{a}$, and $\overline{g}$ is the lifted metric
to the cover (of the same scalar curvature as that of $g$, and twice the 
area), the conformal class $[\overline{g}]$ has 
fundamental domain with interior angles matched in pairs. The 
extremal of the area of these under deformations that preserve these 
pairings is the same as the one described above. Thus, $[\overline{g}]$ is the 
critical class of the restricted Willmore functional $\mc{W}_{2a}$ 
on $\tilde{\Sigma}^{k-1}$, and there is a 2-to-1 immersion 
$\pi: (\tilde{\Sigma}^{k-1},\overline{g}) \rightarrow
f_{g}(\Sigma)\hookrightarrow (\mb{S}^{\tn},\tg)$. If this is the case, 
the critical value is given by the stated expression, by definition.
\qed
\smallskip

By Theorems \ref{th5} and \ref{th7}, we can deduce the existence of 
a preferred scale $a$, namely that given by the area of a minimal 
critical $f_g$ of the functional (\ref{ssss}), in which case 
$\mc{W}_{f_g}(\Sigma)= 4\, {\rm area}_g(\Sigma)=\mc{W}(\Sigma,[g])$. In the 
oriented case, this preferred scale is given by the area 
${\rm area}_{g_{k,1}}(\xi_{k,1})
\in [2\pi^2, 8 \pi)$ of the equilateral Lawson surface.
\smallskip
 
\section{Proof of the main theorem} \label{s3}
In our proof of Theorem \ref{th1}, we parallel the argument we used to study 
strongly extremal metrics, and classes with minimal energy 
\cite{sim5}, \cite{silu1,silu2}. We fix the scale of the metrics by the
preferred scale of a critical isometric embedding across all classes of 
area given by this scale, and for any class represented by a metric of this 
area, use the 
equal area metric of constant scalar curvature in the class, and the extrinsic 
quantities of its isometric embedding to compute the Willmore invariant of 
the class. The parallelism breaks here by the inherent changing complex
structure of the metrics across the classes.
\smallskip

{\it Proof of Theorem \ref{th1}}. We denote by $(\Sigma, g)$ the standard 
sphere $(\mb{S}^2,g)$, or the Lawson surface $(\xi_{k,1},g_{\xi_{k,1}})$ when 
the genus $k$ of $\Sigma$ is zero, or greater, respectively, and let 
$f_g : (\Sigma, g) \rightarrow (\mb{S}^3,\tg)\subset (\mb{S}^{\tn},\tg)$ be its 
minimal isometric embedding. If $J=J_g$ is the complex structure 
that $g$ induces, we refer to the pair $(J,g)$ as the standard K\"ahler 
structure on $\Sigma$. By its minimality, $f_g$
satisfies the system (\ref{wcc}), and we have
$\mc{W}_{f_g}(\Sigma)=\mc{W}(\Sigma,[g])$. Since the group of 
congruences of $(\xi_{k,1},g_{\xi_{k,1}})$ is 
$\mb{Z}/2(k+1) \times \mb{Z}/4$, the fundamental
domain of $[g_{\xi_{k,1}}]$ is equilateral of equal interior angles, and by 
Theorem \ref{th8}, when $k\geq 1$, $[g]$ is the conformal class of the 
critical embedding $f_g$ of $\mc{W}_a$, of critical value 
$4a=4 \, {\rm area}_g(\Sigma)=\mc{W}_{f_g}(\Sigma)=\mc{W}(\Sigma,[g])$.  

Given any metric $g'$ on $\Sigma$, we let $J'$ be the complex structure on 
$\Sigma$ that $g'$ induces, and consider the path of Riemannian metrics  
$$
[0,1] \ni t \rightarrow g_t = (1-t)g' + tg \, .
$$
If $J_t$ is the complex structure on $\Sigma$ that $g_t$ 
induces, we obtain a continuous path $[0,1]\ni t \rightarrow (J_t, g_t)$ of 
K\"ahler structures on $\Sigma$ of K\"ahler classes 
$\Omega_{g_t}= [ \omega_{g_t}]$, $\omega_{g_t}=g_t(J_t \, \cdot \, , \, 
\cdot \,)$, that starts at $(J',g')$ and ends at the standard $(J,g)$. For
$a_t:={\rm area}_{g_t}(\Sigma)$, the solution to the extremal problem of 
Calabi \cite{ca} on $(\Sigma,J_t, g_t)$ produces a minimizer of the functional  
$\mc{S}^2_{a_t,[g_t]}$ in (\ref{ss}) whose
K\"ahler form represents $\Omega_{g_t}$. Notice that this functional is not
conformally invariant, but its minimum is achieved at the end point of a
path of area preserving conformal deformations of $g_t$. We obtain a path
$[0,1]\ni t \rightarrow (J_t,g'_t)$ of Einstein K\"ahler metrics
such that ${\rm area}_{g'_t}(\Sigma)={\rm area}_{g_t}(\Sigma)=a_t$, and 
$g'_t \in [g_t]$. The continuity 
of $t\rightarrow J_t$ reflects the fact that the conformal structures on 
$\Sigma$ are continuous functions of the parameters of its Teichm\"uller space.
By Theorem \ref{th4}, we have that $\mc{W}_{f_{g'_t}}(\Sigma)=
\mc{W}_{f_{g_t}}(\Sigma)$, thus,
$\mc{W}_{f_{g'_0}}(\Sigma)=\mc{W}_{f_{g'}}(\Sigma)$, while
$\mc{W}_{f_{g'_1}}(\Sigma)=\mc{W}_{f_{g}}(\Sigma)$.
By Theorem \ref{th7}, the embedding $f_{g'_t}\mid_{t=1}$ is minimal of 
codimension $1$, and so when the genus $k$ is either $0$, or $1$, we have 
that $(J_1, g'_1)=(J,g)$,
and $f_{g'_1}=f_g$, respectively. 

If there are $t$s such that $\mc{W}_{f_{g'_t}}(\Sigma)>  
\mc{W}(\Sigma,[g_t])$, we can change continuously the path $g'_t$ by 
homothetics to a path of Einstein metrics $\bar{g}'_t$ such that
$\mc{W}_{f_{\bar{g}'_t}}(\Sigma)=\mc{W}(\Sigma,[g_t])$. If 
$\mc{W}(\Sigma,[g_t]) \geq 
\mc{W}(\Sigma,[g])$ for all $t$, we may apply  
Theorems \ref{th4} and \ref{th7} to the metrics along this path, if necessary, 
to conformally dilate them to Einstein metrics $g_t^n$ and  
$g_t^e$, respectively, so that ${\rm area}_{g_t^n}(\Sigma)={\rm area}_g(
\Sigma)$, and if
$$
f_{g_t^n}: (\Sigma,g_t^n)\hookrightarrow (\mb{S}^{\tn},
\tilde{g})
$$
is the corresponding family of Nash isometric embeddings, we have that 
\begin{equation} \label{csc}
s_{g_t^n}=\frac{4\pi \chi}{{\rm area}_g(\Sigma)} =
2+ \| H_{f_{g_t^n}}\|^2 - \| \alpha_{f_{g_t^n}}\|^2 \, ,  
\end{equation}
$g_t^e=(1+\| H_{f_{g_t^n}}\|^2/4)g^n_t$, and the corresponding family of 
isometric embeddings
$$
f_{g_t^e}: (\Sigma,(1+\| H_{f_{g_t^n}}\|^2/4)g_t^n )
\rightarrow (\mb{S}^{2+p},\tilde{g}) \hookrightarrow (\mb{S}^{\tn},\tg)
$$
is minimal for all $t$. The embedding 
$f_{g_t^n}\mid_{t=1}$ is minimal of codimension
$1$, and so when $k$ is either $0$, or $1$, we have that $g_1^n=g_1^e=g$, 
and $f_{g_1^n}=f_{g_1^e}=f_g$, respectively. Since
$\mc{W}_{f_{g_t^n}}(\Sigma)=
\mc{W}(\Sigma,[g_t])$, we have that
$\mc{W}_{f_{g'}}(\Sigma)\geq \mc{W}_{f_{g_0^n}}(\Sigma)$, while  
$\mc{W}_{f_{g_1^n}}(\Sigma)=\mc{W}_{f_{g}}(\Sigma)$.

There are no $t$s such that $\mc{W}_{f_{\bar{g}'_t}}(\Sigma)=
\mc{W}(\Sigma,[g_t]) < \mc{W}(\Sigma,[g])$, and so the construction of 
$g_t^n$ and $g_t^e$ above is always possible. When $k=0$, this follows by
the estimate of Li and Yau on the Willmore energy of any surface without 
self-intersection, which is at least $16\pi=\mc{W}_{f_g}(\Sigma)$ \cite{liya}. 
So assume that $k\geq 1$, 
and that such a $t$ exists. By Theorem \ref{th4}, we may
conformally dilate $\bar{g}'_t$ to an Einstein $h^e_t\in [g_t]$ such that
$f_{h_t^e}$ is minimal, and have 
$\mc{W}_{f_{h^e_t}}(\Sigma)=4{\rm area}_{h_t^e}(\Sigma)=\mc{W}(\Sigma,[g_t])$.
Similarly, we conformally deform $g$, if necessary, to an Einstein $g^e\in [g]$
of the same area, with $f_{g^e}$ minimal, and conformally dilate it to 
an Einstein $\bar{g}^e$ of the same Willmore energy 
such that $a:={\rm area}_{\bar{g}^e}(\Sigma)={\rm area}_{h_t^e}(\Sigma)$. 
Then the functional $\mc{W}_a$ in (\ref{ssss})
has critical points $f_{h_t^e}$ and $f_{\bar{g}^e}$.
This contradicts the assumed relationship between their Willmore enegies, 
since, by Theorem \ref{th8}, the metrics should be both in the conformal 
class of the equilateral Lawson surface.

By (\ref{af}), as the conformal class varies along the path of metrics 
$g_t^n$, the component $T^{\nu}$ of the variational vector field $T$ is 
$L^2$-orthogonal to $H_{f_{g^n_t}}= h_{f_{g^n_t}}\nu_{f_{g^n_t}}$. 
Since all normal directions are conformal, it follows that the functionals 
$\Psi_{f_{g_t^n}}$ and $\Pi_{f_{g_t^n}}$ are stationary when the isometric
embedding is deformed in any normal direction $L^2_{g_t^n}$ orthogonal to 
$\nu_{f_{g_t^n}}$ also. Since the metrics 
$g_t^e=(1+h^2_{f_{g^n_t}}/4)g^n_t$ have their corresponding family of 
isometric embeddings $f_{g_t^e}$ 
minimal for all $t$, the $f_{g_t^e}$s are all critical points 
of $\Psi_{f_{g_t^e}}$ under deformations of the embedding in any 
normal direction, $\nu_{f_{g_t^e}}$ unrestricted included, and it
satisfies the system of equations (\ref{cphi}). 

\begin{enumerate}[label={\rm \alph*)}]
\item Suppose that $k=0$. Given a metric $g'$ on
$\Sigma$, we consider the isometric embeddings $f_{g_t^n}$ and 
$f_{g_t^e}$ constructed
above. The extrinsic functions $t \rightarrow \| \alpha_{f_{g_t^e}}\|^2$, and 
$t \rightarrow \| H_{f_{g^e_t}}\|^2$ both vanish at $t=1$, and by continuity, 
the lower and upper 
estimates in (\ref{es}) hold for a nontrivial open neighborhood of 
$1\in [0,1]$. In this neighborhood, we 
have that  $\| H_{f_{g^e_t}}\|^2=0=\| \alpha_{f_{g^e_t}}\|^2$,
$[g^n_t]=[g_t^e]=[g_t]=[g]$, and   
modulo an isometry, $f_{g^n_t}=f_{g_t^e}$ is the standard 
embedding $f_g: (\mb{S}^2,g)\rightarrow (\mb{S}^{\tn},\tg)$ as a totally
geodesic submanifold \cite[Theorem 9]{scs}. If this neighborhood were not 
equal to the entire interval $[0,1]$, there would be a $\bar{t} \in [0,1)$ at 
which the function $\alpha_{f_{g_{\bar{t}}^e}}$
would satisfy (\ref{esu2}), where $q$ is the codimension of the
embedding $f_{g_{\bar{t}}^e}$. In that case, the   
path of constant functions $[0,1]\ni t \rightarrow \| \alpha_{f_{
g^e_t}}\|^2$, which is identically equal to $0$ nearby $1$, would change
from above $1$ to $0$ somewhere on the interval, contradicting
its continuity. Thus, no such $\bar{t}$ can exist, and  
$[g^e_t]=[g_t^n]=[g_t]=[g]$ and 
$f_{g_t^e}=f_{g_t^n}=f_g$ for all $t\in [0,1]$. Thus, $[g']=[g]$, and
since any of the $f_{g_t^e}$ realizes the  
class invariant $\mc{W}(\Sigma,[g])$, we have that
$$
\mc{W}_{f_{g'}}(\Sigma) \geq \mc{W}_{f_{g_0^e}}(\Sigma)=
\mc{W}_{f_g}(\Sigma)= 4\, \mu_g(\Sigma)= 16 \pi \, .
$$
If the equality on the left holds, 
$f_{g'}$ must have been a conformal deformation of the minimal 
$f_{g_0^n}=f_{g_0^e}=f_g$, of possibly smaller area. 
\item Suppose now that $k\geq 1$. Given a metric $g'$ on $\Sigma$, we 
consider the paths of isometric embeddings $t\rightarrow f_{g_t^n}$ and
$t\rightarrow f_{g_t^e}$ constructed above. If $h_t$ is the mean
curvature function of $f_{g_t^n}$, we have that  
$$
\mc{W}_{f_{g_t^n}}(\Sigma)= 4\, \mu_{g_t^n}(\Sigma)+ 
\Psi_{f_{g_t^n}}(\Sigma) = (4 + h_t^2) \mu_{g_t^n}(\Sigma) \geq 
4\, \mu_g(\Sigma) =
\mc{W}_{f_g}(\Sigma) \, ,
$$
and the lower bound is achieved if, and only if, the isometric embedding
$f_{g_t^n}$ is minimal, in which case it coincides with $f_{g_t^e}$ up to 
isometries. Notice that $\mu_{g_t^e}(\Sigma)= 
(1 + h_t^2/4)\mu_{g_t^n}(\Sigma)$, and 
we see explicitly that the expression on the left of the inequality above 
coincides with $\mc{W}_{f_{g_t^e}}(\Sigma)=\mc{W}(\Sigma,[g_t])$.
Since
$f_{g_0^n}$ and $f_{g_1^n}$ are minimizers in their corresponding classes, we 
conclude that  
$\mc{W}_{f_{g'}}(\Sigma)\geq \mc{W}_{f_{g_0^e}}(\Sigma)$, and 
$\mc{W}_{f_{g}}(\Sigma)=\mc{W}_{f_{g_1^e}}(\Sigma)$, respectively. Thus,
$\mc{W}_{f_{g'}}(\Sigma) \geq 
\mc{W}_{f_{g}}(\Sigma)$, as desired. 

If the equality
$$
\mc{W}_{f_{g'}}(\Sigma) = 4 \mu_g (\Sigma)=\mc{W}_{f_g}(\Sigma)
$$
where to hold, then we would have $\mc{W}_{f_{g'}}(\Sigma)=\mc{W}_{f_{g_0^n}}
(\Sigma)=\mc{W}_{f_{g_0^e}}(\Sigma)=
\mc{W}_{f_{g}}(\Sigma)$, and $f_{g_0^n}$ and $f_{g_0^e}$ would both be
minimal, and equal to each other (up to isometries, of course). If
$a:=\mu_{g_0^e}(\Sigma) = \mu_{g}(\Sigma)$, $f_{g_0^e}$ is an absolute 
minimum of the functional $\mc{W}_a$ in (\ref{ssss}), which, by 
Theorem \ref{th8}, implies that
the conformal class $[g_0^e]=[g_0=g']$ must have fundamental domain 
that is equilateral of equal interior angles. Thus, the metrics 
$g'$ and $g$ must be conformal to each other, and their Nash 
isometric embeddings $f_{g'}$ and $f_g$ differ by at most a
conformal transformation, possibly with ${\rm area}_{g'}(\Sigma) \leq 
{\rm area}_g(\Sigma)$.  
\qed
\end{enumerate}
\smallskip

Corollary \ref{co2} is stated for emphasis, as our argument yields a 
new proof of this classic result. The low end of the gap theorem distinguishes 
$\mb{S}^2$ as a complex manifold with a unique complex structure, while it 
precludes any almost complex structure on $\mb{S}^6$ from being integrable,  
and making this sphere diffeomorphic to a complex 
manifold \cite[Theorem 1]{scs}. 
\smallskip

{\it Proof of Corollary \ref{co2}}. Any metric in the single conformal class 
induces the same complex structure on the (real) surface.   
\qed
\smallskip

The upper end of the gap theorem may be used to characterize the 
minimizer of $\mc{W}$ above in the case when $k=1$, instead of the argument 
presented that works for all $k\geq 1$. This leads to a proof of 
Corollary \ref{co3}, a conjecture of Lawson \cite{law2} from the late 
1960s, recently proven by Brendle \cite{bren} as an application of the maximum 
principle to an elliptic equation that involves $\| \alpha \|$, and which uses 
Simons' identity \cite[Theorem 5.3.1, Lemma 5.3.1]{si} in a crucial manner.  
That our approach could lead into this direction also 
was presaged by the result of Urbano cited earlier, when we outlined Cod\'a 
and Neves' work in \S 1. 
\smallskip

{\it Proof of Corollary \ref{co3}}. We let $g'$ be the intrinsic metric 
on $\Sigma$, and view its minimal embedding as 
$f_{g'}: (\Sigma, g') \rightarrow (\mb{S}^3,\tg) \subset (\mb{S}^{\tn},\tg)$.
There exists a path $[0,1]\ni t \rightarrow g_t$ of area
preserving conformal deformations of $g'$ to an scalar flat metric at $t=1$
that minimizes (\ref{ss}), 
and by Theorem \ref{th7}, the corresponding family $f_{g_t}: (\Sigma,g_t)
\rightarrow (\mb{S}^{\tn},\tg)$ of Nash isometric embeddings is minimal  
and of codimension $1$ for all $t$. By (\ref{sce}), the path 
$t \rightarrow \| \alpha \|^2_{f_{g_t}}$ takes the value $2$ at $t=1$, and
satisfies (\ref{esu1}) for $q=1$. By the gap theorem, up to an isometry of
$\mb{S}^3$,  
$(\Sigma, g_1)$ is the Clifford torus $(\xi_{1,1},g_{\xi_{1,1}})$,  
and $f_{g_1}$ is its standard embedding $f_{g_{\xi_{1,1}}}$.
Thus, $g'\in [g_{\xi_{1,1}}]$,
${\rm area}_{g'}(\Sigma)= {\rm area}_{g_{\xi_{1,1}}}(\xi_{1,1})=
2\pi^2$, and $f_{g'}$ is a minimal conformal deformation of $f_{g_{\xi_{1,1}}}$
in $\mb{S}^3$. 

We reverse the time parameter, and look at the area preserving deformation
$f_{g_t}$ of $f_{g_{\xi_{1,1}}}$ into $f_{g'}$. If we let $u=u(t)$ be the path 
of scalar functions that defines the conformal factor $e^{2u}$ in 
the relation $g_t=e^{2u}g_{\xi_{1,1}}$, by (\ref{eq1}), and (\ref{eq4}), 
we have that
$$
\begin{array}{rcl}
2 & = & e^{-2u}(2+ 2 \Delta^{g_{\xi_{1,1}}} u ) \, , \\
\| \alpha\|^2_{f_{g_t}} & = & 2 e^{-2u} \, , 
\end{array}
$$
respectively. 
In particular, $\| \alpha \|_{f_{g_t}}$ is never zero, and by  Simons' identity
\cite[Theorem 5.3.1]{si}, it follows that
$$
-\Delta^{g_t} \| \alpha \|_{f_{g_t}} -
\frac{ \| \nabla^{g_t} \| \alpha \|_{f_{g_t}}\|^2}{\| \alpha \|_{f_{g_t}}}+
(\| \alpha \|^2_{f_{g_t}}-2)\| \alpha \|_{f_{g_t}}= 0 \, .
$$
We may rewrite this equation in terms of $u$, and $\Delta^{g_{\xi_{1,1}}} u$, 
and obtain that 
$$
\sqrt{2}e^{-u}\left( (e^{-2u}\Delta^{g_{\xi_{1,1}}}u-e^{-2u}\| \nabla^{g_{\xi_{
1,1}}} u\|^2_{g_{\xi_{1,1}}})
e^{-2u}\| \nabla^{g_{\xi_{1,1}}} u\|^2_{g_{\xi_{1,1}}}+
(2e^{-2u}-2)\right)= 0\, ,
$$
which by using the first of the two identities of earlier, 
after a minor simplification, and a cancellation of a nonzero factor,
yields that
$$
\Delta^{g_{\xi_{1,1}}}u+2 \| \nabla^{g_{\xi_{1,1}}}
 u\|^2_{g_{\xi_{1,1}}} = 0 \, .
$$
By integration relative to the measure $d\mu_{g_{\xi_{1,1}}}$, we conclude that
$\nabla^{g_{\xi_{1,1}}}u(t)$ vanishes, which in addition to the fact that
$\nabla u^{\nu}$ vanishes also implies that $u(t)$ is a $t$ dependent constant. 
By the area preserving condition, this 
constant must be zero.
Thus, $u(t)\equiv 0$, and up to an isometry of the background 
$(\mb{S}^3,\tg)$, $g'=g_{\xi_{1,1}}$. 
\qed    
 
\section{The Willmore functional on embeddings of nonoriented surfaces}
\label{s4}
We analyze now the functional (\ref{will}) when the surface
$\Sigma^{k}$ is nonoriented. Although we describe geometrically the    
critical class, and how to find its canonical representative to compute the
critical value abstractly, we are unable to fix the scale explicitly except in 
the cases of topological genus $k=1$ and $k=2$, respectively.

\begin{theorem} \label{th9}
If $\Sigma=\Sigma^k$ is nonoriented, 
and $\mc{W}$ is the Willmore functional {\rm (\ref{will})}, then:
\begin{enumerate}[label={\rm \alph*)}] 
\item If $k=1$, for any Riemannian metric $g'$ on $\Sigma$, we have that  
$$
\mc{W}_{f_{g'}}(\Sigma)\geq \mc{W}_{f_g}(\mb{P}^2(\mb{R}))=24\pi\, ,
$$
where $f_g : (\mb{P}^2(\mb{R}),g) \rightarrow (\mb{S}^4,\tg)$ is the 
minimal isometric embedding of the Einstein $(\mb{P}^2(\mb{R}),g)$ of scalar 
curvature $2/3$, and $g'\in [g]$. The lower bound is achieved by $f_{g'}$ if,
and only if, $f_{g'}$ is conformally equivalent to $f_g$ in 
$(\mb{S}^{\tn},\tg)$, and ${\rm area}_{g'}(\Sigma)\leq 6\pi$.   
\item If $k\geq 2$, and if $g$ is a representative of a class 
$[g]$ with minimal isometric embedding $f_g$, and the lift $\overline{g}$ of 
$g$ to the oriented {\rm 2}-to-{\rm 1} cover 
$\tilde{\Sigma}^{k-1}$ of $\Sigma$ is the conformal class of 
equilateral fundamental domain,  
for any Riemannian metric $g'$ on $\Sigma$, we have that  
$$
\mc{W}_{f_{g'}}(\Sigma) \geq 
\mc{W}_{f_{g}}(\Sigma) = 4\, \mu_{g}(\Sigma) \, ,
$$
and the lower bound is achieved by $f_{g'}$ if, and only if,   
$f_{g'}$ is conformally equivalent to $f_g$ in the ambient space 
background $(\mb{S}^{\tn},\tg)$ and $\mu_{g'}(\Sigma) \leq 
\mu_g(\Sigma)$. If $k=2$, the absolute minimizer is the bipolar Lawson
surface $(\tilde{\tau}_{3,1},g)$ minimally embedded into $(\mb{S}^4,\tg)$,  
of ${\rm area}_g(\tilde{\tau}_{3,1})=6\pi E(2\sqrt{2}/3)$, $E(p)$ the complete 
elliptic integral of parameter $p$.
\end{enumerate}
\end{theorem}

{\it Proof}. a) Over the sphere $\mb{S}^2(\sqrt{3/2})$, the
mapping 
$$
x=[x_0:x_1:x_2] \mapsto f(x)=\frac{1}{\sqrt{3}}
 (2x_0x_1, x_0^2-x_1^2,2x_0 x_2, 2x_1 x_2,
(x_0^2+x_1^2- 2x_2^2)/\sqrt{3}) 
$$
defines a minimal embedding of 
$\mb{P}^2(\mb{R})$ into $(\mb{S}^4,\tg)$, with image of scalar curvature 
curvature $2/3$ \cite[Theorem 1, Remark 2]{si3}. If $g$ denotes the Einstein
metric of scalar curvature $2/3$ on the sphere and projective plane, 
we obtain a commutative diagram  
$$
\begin{array}{rcccc}
\mb{Z}/2 & \hookrightarrow & (\mb{S}^{2}(\sqrt{3/2}),g)  & & \\
& & \pi \downarrow \phantom{m} & \searrow & \\
& & (\mb{P}^2(\mb{R}),g) & \stackrel{f_g}{\rightarrow} & 
(\mb{S}^4,\tg) \, ,  
\end{array}
$$
where $\pi$ is a 2-to-1 covering Riemannian 
submersion. We have that ${\rm area}_g (\mb{P}^2(\mb{R}))=6\pi$, and   
$\mc{W}_{f_{g}}(\Sigma)=24\pi$.

Since the cover carries only one conformal class of metrics, there is
only one conformal class of metrics on $\Sigma$, that given by
$[g]$. By Theorems \ref{th5} and \ref{th7}, $f_g$ is the only minimal 
embedding of an Einstein metric in the class of area $6\pi$, and 
$f_g$ realizes $\mc{W}(\Sigma,[g])$. 

If $f_{g'}$ is any other isometric 
embedding, $g'\in [g]$, and we have that $\mc{W}_{f_{g'}}(\Sigma) \geq 
\mc{W}_{f_{g}}(\Sigma)=24\pi$. If the equality is achieved, 
by Theorem \ref{th4}, 
$f_{g'}$ is a conformal deformation of $f_g: 
(\mb{P}^2(\mb{R}),g) \rightarrow (\mb{S}^4,\tg) \hookrightarrow (\mb{S}^{\tn},
\tg)$ of ${\rm area}_{g'}(\Sigma)\leq 6\pi$.  
\smallskip

b) If $k=2$, $\Sigma$ is the Klein bottle. The manifold $\tau_{3,1}$ 
of Lawson is in the conformal class of the square torus $\mb{R}^2/\Gamma$ of 
lattice $\Gamma$ generated by $(\pi,\pi)$ and $(\pi,-\pi)$, that can be
explicitly described by a doubly periodic immersion into $\mb{S}^3$ 
\cite[(7.1)]{la2}. Its associated bipolar surface 
$\Sigma=\tilde{\tau}_{3,1}$ is a Klein bottle 
minimally embedded in $(\mb{S}^4,\tg)\hookrightarrow (\mb{S}^5,\tg)$, 
whose parametrization can be explicitly described also
\cite[\S3]{lap}. If $g$ is the intrinsic metric on $\tilde{\tau}_{3,1}$, 
we have an isometric minimal embedding
$f_g: (\tilde{\tau}_{3,1},g) \rightarrow (\mb{S}^4,\tg)
\hookrightarrow (\mb{S}^5, \tg)$, and an explicit calculation shows that
$a:={\rm area}_g(\Sigma)= 6\pi E(2\sqrt{2}/3)$. 
If we lift $g$ to a metric $\overline{g}$ on the square torus 
$\mb{R}^2/\Gamma$, 
we have a 2-to-1 covering Riemannian submersion $\pi: (\mb{R}^2/\Gamma,
\overline{g}) \rightarrow (\tilde{\tau}_{3,1},g)$, and the area of the cover 
is $2a$.
This metric is not Einstein, but by Theorem \ref{th5}, its equal area 
Einstein deformation has the same Willmore energy. 
By Theorem \ref{th7}, $a$ is the only scale for
which there are metrics in $[g]$ that are minimally embedded, and by 
Theorem \ref{th8}, $f_g$ realizes $\mc{W}(\Sigma,[g])$.   

If $f_{g'}$ is any other isometric embedding, by Theorem \ref{th8} applied
to the restricted Willmore functional $\mc{W}_a$, we conclude that 
$\mc{W}_{f_{g'}}(\Sigma) \geq 
\mc{W}_{f_{g}}(\Sigma)=24\pi E(2\sqrt{2}/3)=4a$. If the equality is achieved, 
by Theorem \ref{th4}, 
$f_{g'}$ is a conformal deformation of $f_g: 
(\tilde{\tau}_{3,1},g) \rightarrow (\mb{S}^4,\tg) \hookrightarrow 
(\mb{S}^{\tn}, \tg)$ of ${\rm area}_{g'}(\Sigma)\leq a$.  

For $k>2$ the argument above generalizes verbatim under the stated hypothesis,
the exception being the knowledge of the explicit value of the scale 
${\rm area}_g(\Sigma)$.
\qed 
\smallskip

\begin{remark}
If $k\geq 3$, the nonoriented minimal surface 
$\eta_{k-1,1}=\iota_{\eta_{k-1,1}}(\Sigma) \hookrightarrow \mb{S}^3$ of 
Lawson has topological genus $k$, immersed but not embedded. Its 
group of congruences is sufficiently large so that the conformal class
of the associated 
bipolar surface $\tilde{\eta}_{k-1,1} \hookrightarrow (\mb{S}^5,\tg)$,
which is nonsingular and minimally immersed, satisfies the conditions of
Theorem \ref{th8}, hence, a plausible candidate for the minimizer of
the Willmore functional on $\Sigma^k$. But we cannot see why the immersion
would be an embedding, which would be needed to draw that conclusion on
the basis of the said theorem.     
\end{remark}

\begin{remark}
In terms of their topological genus, the first two 
minimums of the Willmore functional on oriented and nonoriented surfaces
are interlaced as 
$$
16 \pi < 24\pi < 8 \pi^2 < 24\pi E(2\sqrt{2}/3)\, ,   
$$
while those in the oriented case for all the remaining $k$s,  
$4 {\rm area}_{\xi_{k,1}}(\xi_{k,1})$, forms an increasing sequence
bounded above by $32\pi$. 
It is tempting to conjecture that if $m_{k}$ is the minimum of the
Willmore functional of a nonoriented surface of topological genus $k$,
then 
$$
4 {\rm area}_{\xi_{k-1,1}}(\xi_{k-1,1}) < m_k <
4 {\rm area}_{\xi_{k,1}}(\xi_{k,1}) \, , \quad k=2,3, \ldots 
$$
We lack sufficient information on the scales $m_k$ to ellucidate this. 
\end{remark}

\end{document}